\input amstex
\input pictex
\documentstyle{amsppt}
\PSAMSFonts
\nologo

\def\ZZ{{\Bbb Z}}

\def\RR{{\Bbb R}}

\def\boxx{\unskip \nopagebreak \hfill $\square$}

\rightheadtext{Metrics and embeddings}

\topmatter
\title
Metrics and embeddings of generalizations\\
of Thompson's group $F$
\endtitle
\author
J. Burillo, S. Cleary and M.I. Stein
\endauthor
\abstract
The distance from the origin in the word metric for generalizations $F(p)$ of
Thompson's group
$F$ is quasi-isometric to the number of carets in the reduced rooted tree
diagrams representing the elements of $F(p)$.  This interpretation of the
metric is used to prove that every $F(p)$ admits a quasi-isometric embedding
into every $F(q)$, and also to study the
behavior of the shift maps under these embeddings.
\endabstract
\address
Dept. of Mathematics, Tufts University, Medford, MA 02155, U.S.A.
\endaddress
\email
jburillo\@emerald.tufts.edu
\endemail 
\address
Dept. of Mathematics, California State University, Fresno, CA 93710, U.S.A.
\endaddress
\email
sean\_cleary\@csufresno.edu
\endemail
\address
Dept. of Mathematics, Trinity College, Hartford, CT 06106, U.S.A.
\endaddress
\email
mstein\@mail.trincoll.edu
\endemail
\endtopmatter

\document

When obtaining the first example of a finitely presented 
infinite simple group in \cite{12}, R. Thompson introduced  groups
$F$,
$T$ and
$V$, which have attracted considerable interest since then and
have become known as Thompson's groups $F$,
$T$ and
$V$.
Higman  \cite{10} generalized the group $V$ to a whole family of
groups, a generalization that was extended to the subgroups $F$ and
$T$ by Brown \cite{3}. The groups $F(p)$ of this article correspond
to the groups $F_{p,\infty}$ in \cite{3}. These families of groups have also 
been extensively studied by Brin and Guzm\'an \cite{2}.

Thompson's group $F$ has proven to be
 the most interesting group of
all and has emerged in a variety of settings.  It can be regarded
as the group of piecewise-linear, orientation-preserving
homeomorphisms of the interval
$[0,1]$ which have breakpoints only at dyadic points and
whose slopes, where defined, are powers of two.  In 1984 Brown and
Geoghegan \cite{4} found it to be the first example  of a finitely
presented infinite-dimensional torsion-free
$FP_{\infty}$ group. This fact has been extended to all $F(p)$ in
\cite{3}, and studied by Stein \cite{11}, where generalizations to more
general  groups of homeomorphisms with general rational breakpoint
sets are considered.  Cleary \cite{7} has studied these properties for
groups of piecewise-linear homeomorphisms with irrational breakpoints 
and slopes.

In \cite{5} Burillo showed the relationship between the word metric of
Thompson's group $F$ and an estimate of the distance derived from the
unique normal form of the elements. This algebraic estimate is 
quasi-isometric to
the word metric and was used to prove that some subgroups are nondistorted in 
$F$. In this paper we find a geometric estimate of the word metric
in terms of rooted tree diagrams (see \cite{6} and 
\cite{8}), show that this is quasi-isometric to the word metric,
and use this
interpretation to prove that the embeddings of $F(p)$ into $F(q)$ outlined 
in \cite{2} are
quasi-isometric. The interpretations of these embeddings in terms 
of the tree diagrams
also yield insights into the behavior of the shift maps under the embeddings.

After a statement of some results about the groups $F(p)$ in section 1 that 
will be needed later in the paper, including a brief description of the
rooted tree diagrams, section 2 contains the interpretation of
the word metric in terms of the normal form and of the number of carets in
the tree diagrams. In section 3 the natural embedding of $F(p^k)$ into $F(p)$
is studied and proved to be nondistorted, and in section 4 the
embeddings of any $F(p)$ in any $F(q)$ obtained in \cite{2} are proved to 
be quasi-isometric.
The last section is dedicated to study the
behavior of the shift maps of $F(p)$ under these embeddings.

\heading
1. The groups $F(p)$
\endheading

The generalizations of Thompson's group which are the subject of this
paper are the groups $F(p)=F(1,\ZZ[\frac 1p],\left<p\right>)$, the groups of
piecewise-linear, orientation-preserving homeomorphisms of the interval
$[0,1]$ which have breakpoints only in $\ZZ[\frac 1p]$, and such that the slopes,
where defined, are powers of $p$. For $p=2$ this group is
the well-known Thompson's group $F$, and the groups $F(p)$ are natural
generalizations of $F$, and share many of its
interesting properties. In this section we outline some of these
properties that will be used later in the paper. For a very readable
introduction to $F$ see \cite{6}, and for generalizations
to $F(p)$ see \cite{3} and \cite{11}.

The group $F(p)$ admits the following infinite presentation:
$$
\Cal P_p=\left<x_i,\,i\ge 0\,|\,x_i^{-1} x_j
x_i=x_{j+p-1},\text{ for }i<j\right>
$$
where the maps $x_i\in F(p)$, for $0\le i\le p-2$, can be represented by the
homeomorphisms of the unit interval:
$$
x_i(t)=\cases
t&\text{if }0\le t\le\frac ip,\\
pt-\frac{i(p-1)}p&\text{if }\frac ip\le t\le\frac{(p-1)(i+1)}{p^2},\\
t+\frac{(p-1)(p-i-1)}{p^2}&\text{if }
\frac{(p-1)(i+1)}{p^2}\le t\le\frac{i+1}p,\\
\frac{t+p-1}p&\text{if }\frac{i+1}p\le t\le 1;
\endcases
$$

For $x_i$ with $i\ge p-1$, we let $j=\left\lfloor\frac i{p-1}\right\rfloor$,
and let $k=i-j(p-1)$. Then,
$x_i\in F(p)$ is the identity except in the interval $\left[1-\frac
1{p^j},1 \right]$, where the graph is a scaled-down copy of the graph of
$x_k$. The compositions are taken on the right; that is, the element
$x_ix_j\in F(p)$ corresponds with the composition $x_j\circ x_i$ as
maps in $[0,1]$.

The groups $F(p)$ admit a {\it shift map\/} $\phi$, which sends every
generator $x_i$ to $x_{i+1}$. The shift map satisfies
$x_0^{-1}\phi(x)x_0=\phi^p(x)$ for all $x\in F(p)$, so it is a conjugacy
idempotent. The relations between the shift maps of the different $F(p)$
are studied in section 5.

The infinite presentation $\Cal P_p$ is useful for its symmetry and simplicity, but
to discuss the word metric we need to consider a finite presentation. It is clear that
$x_0,x_1,\ldots,x_{p-1}$ generate $F(p)$, and it is possible to write a
presentation for $F(p)$ with these $p$ generators and $p(p-1)$
relators (see \cite{9} and \cite{11}). 
For $p=2$ this is the standard presentation
for Thompson's group $F$:
$$
\left<x_0,\,x_1\,|\,[x_0x_1^{-1},x_2],\,[x_0x_1^{-1},x_3]\right>.
$$

In the following, when we refer to the word metric, or the length of
an element, we will always mean with respect to these finite
presentations.

From the relators $x_i x_j x_i^{-1}=x_{j+p-1}$ we
can see that every element in $F(p)$ admits an
expression of the form
$$
x_{i_1}^{r_1} x_{i_2}^{r_2}\ldots x_{i_n}^{r_n}
x_{j_m}^{-s_m}\ldots x_{j_2}^{-s_2} x_{j_1}^{-s_1}
$$
where $i_1<i_2<\ldots<i_n\ne j_m>\ldots>j_2>j_1$. This expression is
unique if we require one additional condition: if both $x_i$ and
$x_i^{-1}$ appear, then one of the generators
$$
x_{i+1},\,x_{i+1}^{-1},\,x_{i+2},\,x_{i+2}^{-1},\ldots,
x_{i+p-1},\,x_{i+p-1}^{-1}
$$
must appear as well.  This is required for uniqueness because if none of them appears,
there is a subword of the type $x_i\phi^p(x)x_i^{-1}$ which can be
replaced by $\phi(x)$. This unique expression of an element will be called
its {\it normal form\/}. The proof of the uniqueness of the normal
form in $F(2)$ in \cite{4}  extends easily to every $F(p)$.
Given an element $x$, its normal form is the shortest word in the
infinite generating set of $\Cal P_p$ representing it.

The groups $F(p)$ can also be realized as groups of homeomorphisms of the
real line $\RR$. The maps
$$
f_i(t)=\cases
t&\text{if }t\le i,\\
pt+i(1-p)&\text{if }i\le t\le i+1,\\
t+p-1&\text{if }i+1\le t
\endcases
$$
generate a group of piecewise-linear homeomorphisms of $\RR$ which is
isomorphic to $F(p)$. As before, compositions are taken on the right.
We will use both of these geometric representations for $F(p)$ to deduce
different properties of these groups.

Another interpretation for $F(p)$ is based on maps of
rooted trees. This interpretation was studied by Higman in \cite{10} 
and Bieri--Strebel in \cite{1}, and it is described with great detail for
the case of $F(2)$ in \cite{6}. A rooted $p$-tree is a tree with a
distinguished vertex (the {\it root\/}) which has $p$ edges, and any
other vertex has either valence 1 ({\it leaves\/}) or valence $p+1$ ({\it
nodes\/}). We think of a rooted $p$-tree as a {\it descending\/} tree,
with the root on top, and different levels of vertices, with the root
being the sole vertex of level 0. Every vertex different from the root is
connected by an edge to a vertex in the previous level, and it is either a
leaf, in which case it is not connected to the next lower level, or a node,
which has $p$ {\it children\/}, i.e., it is connected to $p$
vertices in the next lower level. A node, together with its $p$ children, and
the $p$ edges connecting them, form a {\it caret\/}. A caret has, 
from left to right in the obvious way, a left edge, 
several interior edges (none if $p=2$),
and a right edge.

\beginpicture
\setcoordinatesystem point at -180 170
\setplotarea x from -180 to 180, y from -143 to 170
\putrule from -150 60 to -60 60
\putrule from 60 60 to 150 60
\putrule from -150 57 to -150 63
\putrule from -120 57 to -120 63
\putrule from -110 57 to -110 63
\putrule from -100 57 to -100 63
\putrule from -90 57 to -90 63
\putrule from -60 57 to -60 63
\putrule from 150 57 to 150 63
\putrule from 120 57 to 120 63
\putrule from 140 57 to 140 63
\putrule from 130 57 to 130 63
\putrule from 90 57 to 90 63
\putrule from 60 57 to 60 63
\putrectangle corners at -45 -90 and 45 0
\putrule from -15 -93 to -15 -87
\putrule from -5 -93 to -5 -87
\putrule from 5 -93 to 5 -87
\putrule from 15 -93 to 15 -87
\putrule from -48 -60 to -42 -60
\putrule from -48 -30 to -42 -30
\putrule from -48 -20 to -42 -20
\putrule from -48 -10 to -42 -10
\linethickness=2pt
\plot -135 120 -105 150 /
\plot -105 150 -75 120 /
\plot 125 120 95 150 /
\plot 95 150 65 120 /
\plot -105 150 -105 90 /
\plot 95 150 95 120 /
\plot -125 90 -105 120 /
\plot -105 120 -85 90 /
\plot 105 90 125 120 /
\plot 125 120 125 90 /
\plot 125 120 145 90 /
\plot -45 -90 -15 -60 /
\plot -15 -60 -5 -30 /
\plot -5 -30 15 -10 /
\plot 15 -10 45 0 /
\put {0} [t] at -150 54
\put {0} [t] at 60 54
\put {0} [rt] at -48 -93
\put {1} [t] at -60 54
\put {1} [t] at 150 54
\put {1} [t] at 45 -93
\put {1} [r] at -48 0
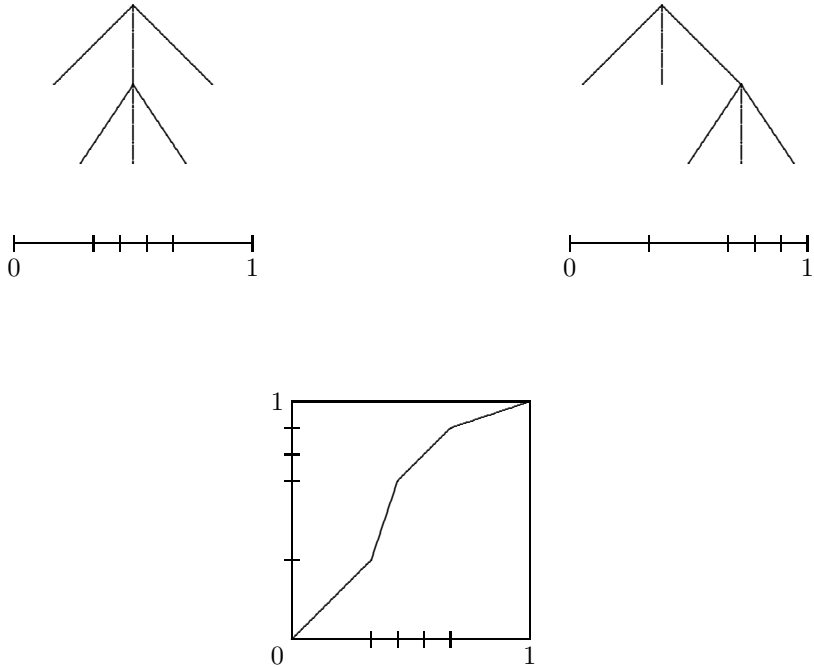
\put {{\bf Fig. 1:} The tree diagram for $x_1\in F(3)$ with its
homeomorphism of $[0,1]$.} [t] at 0 -113
\endpicture

Following \cite{8}, we will say that a caret is a 
{\it left caret\/} when it is the leftmost caret of a level, or, equivalently, 
when its vertex is joined to the root by a path made out completely by left 
edges. We define a {\it right caret\/} in a similar manner, 
and every caret which is 
not left or right is called an {\it interior caret.}

A rooted $p$-tree can be thought of as a graphic representation of a
subdivision of the interval $[0,1]$ into intervals of lengths
$\frac1{p^k}$ for different $k$. A vertex in level $k$ corresponds
to an interval of length $\frac1{p^k}$. If it is a node, the caret
represents the subdivision of the interval in $p$ intervals of length
$\frac1{p^{k+1}}$. The leaves represent the final intervals of the
subdivision, and the order of the intervals in $[0,1]$ induces a total
order on the leaves of the tree. The endpoints of the intervals are always
elements of $\ZZ[\frac1p]$.

\beginpicture
\setcoordinatesystem point at -180 80
\setplotarea x from -180 to 180, y from -115 to 80
\linethickness=2pt
\plot -130 30  -100 60 /
\plot -120 40  -110 30 /
\plot -110 50  -100 40 /
\plot -100 60  -90 50 /
\plot -60 50  -45 60 /
\plot -45 60  -30 50 /
\plot -50 30  -30 50 /
\plot -30 50  -20 40 /
\plot -40 40  -30 30 /
\plot 20 40  30 50 /
\plot 30 50  40 40 /
\plot 30 50  45 60 /
\plot 45 60  60 50 /
\plot 90 50  100 60 /
\plot 100 60  120 40 /
\plot 100 40  110 50 /
\plot -60 -40  -20 0 /
\plot -20 0  -10 -10 /
\plot -30 -10  -20 -20 /
\plot -40 -20  -30 -30 /
\plot -50 -30  -40 -40 /
\plot 10 -10  20 0 /
\plot 20 0  50 -30 /
\plot 20 -20  30 -10 /
\plot 20 -40  40 -20 /
\plot 30 -30  40 -40 /
\setdashes <2pt>
\plot -140 20  -130 30 /
\plot -130 30  -120 20 /
\plot -70 40  -60 50 /
\plot -60 50  -50 40 /
\plot 40 30  60 50 /
\plot 60 50  70 40 /
\plot 50 40  60 30 /
\plot 100 20  120 40 /
\plot 120 40  130 30 /
\plot 110 30  120 20 /
\put {$x_0^2x_1^{-1}$} [t] at -80 20
\put {$x_0$} [t] at 80 20
\put {$x_0^2x_1^{-1}x_0=x_0^3x_2^{-1}$} [t] at 0 -45
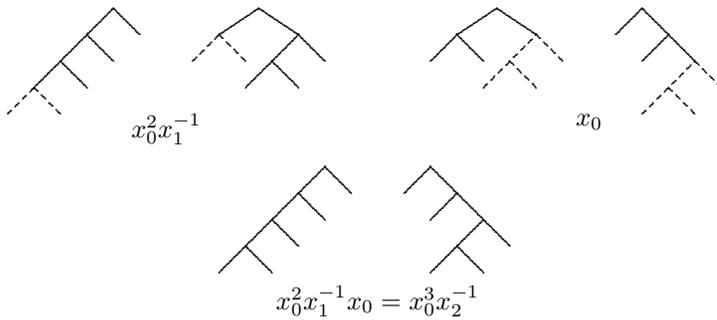
\put {{\bf Fig. 2:} Composing two tree diagrams in $F(2)$, with the 
added carets in dashes.} [t] at 0 -75
\endpicture

An element $x\in F(p)$ can be thought of, then, as a {\it tree
diagram\/}, i.e., an order-preserving map between the leaves of two rooted
$p$-trees with the same number of carets (and thus the same number of leaves). 
The
homeomorphism $x$ of $[0,1]$ represented by this diagram
can be seen using the subdivisions represented by the two trees: the
intervals of the source tree are sent to the ones in the target tree,
preserving the order. The product of two elements $x$ and $y$ can be
seen with the tree diagrams using the following process: add 
carets simultaneously at corresponding leaves to the source  and target trees of $x$, and also
add carets to the pair of trees of $y$,
until the target tree of $x$ and the
source tree of $y$ are equal. Then, a tree diagram for $xy$ has for
source the source of $x$ and for target the target of $y$, with all the
added carets needed to perform the composition.

A diagram is reducible if all the leaves of a caret are sent to all the
leaves of another caret in the target, that is, if these carets
represent
superfluous subdivisions of the corresponding intervals.  In the following, 
we will assume that
our diagrams are reduced, meaning that superfluous carets have been eliminated. For a given
element, the reduced tree diagram representing it is unique, and there is a close
relation between the reduced tree diagram and the normal form of the
element (\cite{6}, \cite{8}).

Reduced $p$-tree diagrams are a powerful and efficient way to represent
elements of $F(p)$, and they will be used several times later in the
paper. A complete detailed description of the tree diagram representation
for $F(2)$ can also be found in \cite{8}, where B. Fordham uses it
to obtain an algorithm to compute the exact length of an element of
$F(2)$ given its normal form.

Yet a different representation of the groups $F(p)$ can be found in
\cite{9}, in the context of diagram groups. This representation is
essentially the same as the tree diagram representation, where a node
with its $p+1$ edges is replaced by a 2-cell with its boundary subdivided
in $p+1$ edges, and the cells are attached to each other along the
edges according to the same pattern represented by the trees.

\heading
2. The word metric of $F(p)$
\endheading

The different interpretations of the groups $F(p)$, both as groups of
homeomorphisms and as groups of maps of rooted trees can be used to
deduce expressions for the word metric. First we 
generalize to all $F(p)$ the estimate of the length of an element developed in \cite{5}.
This gives a quantity which is equivalent to the length, and can be
readily computed from the
normal form of an element.  Here, we mean equivalent in the sense of
quasi-isometry.  We will denote by $|x|_p$
the minimal length of an element of $F(p)$ in the word metric with respect
to  the generators
$x_0,x_1,\ldots,x_{p-1}$.

\proclaim{Theorem 1} Let $x\in F(p)$ have normal form
$$
x_{i_1}^{r_1} x_{i_2}^{r_2}\ldots x_{i_n}^{r_n}
x_{j_m}^{-s_m}\ldots x_{j_2}^{-s_2} x_{j_1}^{-s_1},
$$
and let
$$
D(x)=r_1+r_2+\ldots+r_n+s_1+s_2+\ldots+s_m+i_n+j_m.
$$
Then, we have
$$
\frac{D(x)}{3(p-1)}\le |x|_p\le 3D(x).
$$
\endproclaim

\demo{Proof} The upper bound is straightforward: replacing each $x_i$
in the normal form
by its expression in the generators $x_0,x_1,\ldots,x_{p-1}$, the
length of the word obtained is less than $3D(x)$, and clearly it is an
upper bound for the minimal length.

For the lower bound, observe that since the normal form has the shortest
possible length among words representing $x$ in the generators of $\Cal
P_p$, we can conclude that $|x|_p\ge
r_1+\ldots+r_n+s_1+\ldots+s_m$. Finally, let $w$ be the word in the
generators $x_0,x_1,\ldots,x_{p-1}$ which has minimal length $|x|_p$. To
obtain the unique normal form from $w$, the generators have to be
switched around using the relators of $\Cal P_p$, at the price of
increasing the index of one of them 
by $p-1$ per switch. A given generator in $w$
needs to be switched at most $|x|_p-1$ times, so the highest
possible generator appearing in the normal form has index at most
$$
p-1+(p-1)(|x|_p-1)=(p-1)|x|_p
$$
and then
$$
i_n\le(p-1)|x|_p\qquad\text{and}\qquad j_m\le(p-1)|x|_p.
$$
Combining all these inequalities we obtain the desired lower bound.\boxx
\enddemo

Part of the previous proof is due to V. Guba, who improved on the proof
given in \cite{5}.

The number $D(x)$ given in the previous result is, then, equivalent to
the distance.  This readily computable quasi-metric $D(x)$ can be used in the place of the
genuine word metric to obtain geometric
characterizations for the distance.

\proclaim{Proposition 2} Let $x\in F(2)$ be an element whose normal form is
a {\it positive\/} word, i.e.,
$$
x=x_{i_1}^{r_1} x_{i_2}^{r_2}\ldots x_{i_n}^{r_n}.
$$
Let
$$
N(x)=\max\,\{i_k+r_k+r_{k+1}+\ldots+r_n+1,\text{ for }k=1,2,\ldots,n\}.
$$
Then,
\roster
\item $N(x)$ is equal to the number of carets in either tree of the
reduced $2$-tree diagram for $x$,
\item $N(x)$ is equal to the $y$-coordinate of the last breakpoint of the
graph of $x$ represented as a homeomorphism of $\RR$,
\item $N(x)$ is quasi-equivalent to the distance. In particular,
$$
\frac{D(x)}2\le N(x)\le D(x)+1.
$$
\endroster
\endproclaim

\demo{Proof} The statement that $N(x)$ is equivalent to the distance is
clear from the definition: clearly $N(x)\le D(x)+1$, and also
$$
\gather
N(x)\ge i_1+r_1+\ldots+r_n+1\ge r_1+\ldots+r_n,\text{ and}\\
N(x)\ge i_n+r_n+1\ge i_n,
\endgather
$$
from which the lower bound follows.

Both statements (1) and (2) will be proved using the same induction in
$n$. If $n=1$, then $x=x_i^r$, and we have that the tree diagram of $x$
has exactly $i+r+1$ carets (see Figure 3).
The graph for the homeomorphism representing $x$ has the following
slopes at points with the given $y$-coordinate:

$$
\matrix \format\c&\quad\l\\
1&\text{for }y\in(-\infty,i),\\
2^r&\text{for }y\in(i,i+2),\\
2^{r-k}&\text{for }y\in(i+k+1,i+k+2),\quad k=1,\ldots,r-1,\\
1&\text{for }y\in(i+r+1,\infty).
\endmatrix
$$
In particular, the last breakpoint is $(i+1,i+r+1)$.

\beginpicture
\setcoordinatesystem point at -180 105
\setplotarea x from -180 to 180, y from -50 to 105 
\plot -87 93 -85 95  -65 75 /
\plot -55 65  -35 45 -37 43 /
\plot -63 43 -65 45  -80 30 /
\plot -90 20  -105 5 -103 3 /
\plot 43 73 45 75  65 55 /
\plot 75 45  95 25 93 23 /
\setdots <3pt>
\plot -70 70  -50 50 /
\plot -70 30  -90 10 /
\plot 60 50  80 30 /
\setsolid
\plot -100 80  -90 90 /
\plot -90 90  -70 70 /
\plot -90 70  -80 80 /
\plot -70 30  -50 50 /
\plot -50 50  -40 40 /
\plot -60 40  -50 30 /
\plot -100 0  -90 10 /
\plot -90 10  -80 0 /
\plot 30 60  40 70 /
\plot 40 70  60 50 /
\plot 40 50  50 60 /
\plot 70 20  80 30 /
\plot 80 30  90 20 /
\put {$i+1$ carets} [lb] at -60 70
\put {$r$ carets} [rb] at -85 25
\put {$i+r+1$ carets} [lb] at 70 50
\put {{\bf Fig. 3:} The tree diagram for $x_i^r$.} [t] at 0 -10
\endpicture

For $n>1$, let $x=x_{i_1}^{r_1} x_{i_2}^{r_2}\ldots x_{i_n}^{r_n}$, and
take $y=xx_i^r$. When composing the two tree diagrams of $x$ and $x_i^r$,
we need to add carets to the middle two trees to make them match. We need
to study two cases.

{\it Case 1:\/} Suppose that $i+1\le N(x)$. 

By induction hypothesis, the trees for the reduced tree diagram for
$x$ have $N(x)$ carets,
and we need to add $r$ carets to the target tree of $x$, and $N(x)-i-1$ 
to the source tree for $x_i^r$ to perform the composition. In any case 
the resulting trees have $N(x)+r$ carets, and observe that if 
$i+1\le N(x)$, then $N(y)=N(x)+r$. It is important to realize that
the tree diagram obtained for $y$ is reduced. This fact can be observed 
using the relations between the reduced tree diagram and the exponents
which appear in the unique normal form (see \cite{8}).

\beginpicture
\setcoordinatesystem point at -180 110
\setplotarea x from -180 to 180, y from -60 to 110
\plot -120 20 -105 50 -90 20 /
\put {.} at -100 25
\put {.} at -105 25
\put {.} at -110 25
\plot -65 80 -60 90 -50 70 /
\plot -60 70 -55 80 /
\plot -50 40 -40 50 -30 40 -25 30 /
\plot -35 30 -30 40 /
\plot -25 10 -20 20 -15 10 /
\plot 20 80 25 90 35 70 /
\plot 25 70 30 80 /
\plot 25 20 35 40 45 50 55 40 /
\plot 30 30 35 20 /
\plot 35 40 40 30 /
\plot 15 0 20 10 25 0 /
\plot 90 70 95 80 105 60 /
\plot 95 60 100 70 /
\plot 120 10 125 20 130 10 /
\setdots <2pt>
\plot -50 70 -40 50 /
\plot -65 10 -60 20 /
\plot -25 30 -20 20 /
\plot 35 70 45 50 /
\plot 20 10 25 20 /
\plot 60 30 65 20 /
\plot 105 60 125 20 /
\setdashes <2pt>
\plot -70 0 -65 10 -60 0 /
\plot -60 20 -50 40 -45 30 /
\plot -55 30 -50 20 /
\plot 50 30 55 40 60 30 /
\plot 60 10 65 20 70 10 /
\put {$x$} [t] at -80 -5
\put {$x_i^r$} [t] at 80 -5
\put {{\bf Fig. 4:} The composition of $x$
with $x_i^r$ when $i+1\le N(x)$.} [t] at 0 -30
\endpicture

For the graphs of the corresponding homeomorphisms, again from 
$i+1\le N(x)$, we see that when composing $x$ with $x_i^r$ the 
only interval that gets modified is $[i,i+1]$, which gets stretched 
into $[i,i+r+1]$. The last breakpoint of $y$ is now $N(x)+r$, which as 
before is equal to $N(y)$.

\beginpicture
\setcoordinatesystem point at -180 120
\setplotarea x from -180 to 180, y from -50 to 120
\putrule from -95 10 to -30 10
\putrule from -85 0 to -85 80 
\putrule from 30 10 to 105 10 
\putrule from 40 0 to 40 100 
\plot -87 60 -83 60 /
\plot -87 40 -83 40 /
\plot -87 30 -83 30 /
\plot 38 80 42 80 /
\plot 38 60 42 60 /
\plot 38 30 42 30 /
\linethickness=2pt
\plot -95 0 -85 10 /
\plot -70 30 -65 40 /
\plot -50 60 -30 80 /
\plot 30 0 40 10 /
\plot 55 30 60 40 /
\plot 65 50 70 60 /
\plot 85 80 105 100 /
\setdashes <2pt>
\plot -85 10 -70 30 /
\plot -65 40 -50 60 /
\plot 40 10 55 30 /
\plot 70 60 85 80 /
\setdots <2pt>
\plot 60 40 65 50 /
\plot -83 60 -50 60 /
\plot -83 40 -65 40 /
\plot -83 30 -70 30 /
\plot 42 80 85 80 /
\plot 42 60 70 60 /
\plot 42 30 55 30 /
\put {$x$} [t] at -62 0
\put {$xx_i^r$} [t] at 67 0
\put {$i$} [r] at -90 30
\put {$i+1$} [r] at -90 40
\put {$N(x)$} [r] at -90 60
\put {$i$} [r] at 35 30
\put {$i+r+1$} [r] at 35 60
\put {$N(x)+r$} [r] at 35 80
\put {{\bf Fig. 5:} The graphs of $x$ and $xx_i^r$ when $i+1\le N(x)$.}
[t] at 0 -20
\endpicture

{\it Case 2:\/} If $i+1>N(x)$, then in the composition of the two tree 
diagrams we need to add $i+r+1-N(x)$ carets to the target tree of $x$ to 
match the source of $x_i^r$: 

\beginpicture
\setcoordinatesystem point at -180 110
\setplotarea x from -180 to 180, y from -60 to 110
\plot -100 20 -85 50 -70 20 /
\plot -50 80 -45 90 -40 80 /
\put {.} at -80 25
\put {.} at -85 25
\put {.} at -90 25
\plot -40 60 -35 70 -30 60 /
\plot 15 80 20 90 25 80 /
\plot 25 60 30 70 40 50 /
\plot 30 50 35 60 /
\plot 35 20 45 40 50 30 /
\plot 40 30 45 20 /
\plot 25 0 30 10 35 0 /
\plot 70 70 75 80 85 60 /
\plot 75 60 80 70 /
\plot 100 10 105 20 110 10 /
\setdashes <2pt> 
\plot -35 50 -30 60 -25 50 /
\plot -30 20 -20 40 -15 30 /
\plot -25 30 -20 20 /
\plot -40 0 -35 10 -30 0 /
\setdots <2pt>
\plot -40 80 -35 70 /
\plot -25 50 -20 40 /
\plot -30 20 -35 10 /
\plot 25 80 30 70 /
\plot 40 50 45 40 /
\plot 35 20 30 10 /
\plot 85 60 105 20 /
\put {$x$} [t] at -60 -5
\put {$x_i^r$} [t] at 60 -5
\put {{\bf Fig. 6:} The composition of $x$ with $x_i^r$ when $i+1>N(x)$.}
[t] at 0 -30
\endpicture

The resulting trees have $i+r+1$ carets, and if $i+1>N(x)$, it is clear 
that $N(y)=i+r+1$ as well. For the homeomorphisms, if $i+1>N(x)$, all the
modifications to the graph of $x$ occur above the last breakpoint, so
the last breakpoint for $y$ is the same than for $x_i^r$, which has
$y$-coordinate $i+r+1$.\boxx
\enddemo

\beginpicture
\setcoordinatesystem point at -180 125
\setplotarea x from -180 to 180, y from -50 to 125
\putrule from -90 10 to -30 10 
\putrule from -80 0 to -80 85 
\putrule from 30 10 to 90 10 
\putrule from 40 0 to 40 105 
\putrule from -82 60 to -78 60  
\putrule from -82 50 to -78 50  
\putrule from -82 70 to -78 70 
\putrule from 38 90 to 42 90 
\putrule from 38 60 to 42 60 
\putrule from 38 50 to 42 50 
\plot -90 0 -80 10 /
\plot -65 50 -30 85 /
\plot 30 0 40 10 /
\plot 55 50 65 60 /
\plot 65 60 68.3 70 /
\plot 71.7 80 75 90 /
\plot 75 90 90 105 /
\setdashes <2pt>
\plot -80 10 -65 50 /
\plot 40 10 55 50 /
\setdots <2pt>
\plot 68.3 70 71.7 80 /
\plot -78 50 -65 50 /
\plot -78 60 -55 60 /
\plot -78 70 -45 70 /
\plot 42 90 75 90 /
\plot 42 60 65 60 /
\plot 42 50 55 50 /
\put {$x$} [t] at -60 0
\put {$xx_i^r$} [t] at 60 0
\put {$i$} [r] at -85 60
\put {$i+1$} [r] at -85 70
\put {$N(x)$} [r] at -85 48
\put {$N(x)$} [r] at 35 48
\put {$i$} [r] at 35 60
\put {$i+r+1$} [r] at 35 90
\put {{\bf Fig. 7:} The graphs of $x$ and $xx_i^r$ when $i+1>N(x)$.}
[t] at 0 -20
\endpicture

For the case of a general word, not necessarily positive, there is no 
relation between the $y$-coordinate of the last breakpoint and the
distance. The elements
$$
x_0x_1\ldots x_{k-1}x_k^2x_{k+1}^{-1}x_k^{-1}\ldots x_1^{-1}x_0^{-1}
$$
have all breakpoints in the square $[0,1]\times[0,1]$, whereas their
norm is linear in $k$. But the number of carets in the reduced diagram 
is still equivalent to the norm.

\proclaim{Theorem 3} Let $x\in F(2)$ be an element whose normal form is
$$
x_{i_1}^{r_1} x_{i_2}^{r_2}\ldots x_{i_n}^{r_n}
x_{j_m}^{-s_m}\ldots x_{j_2}^{-s_2} x_{j_1}^{-s_1},
$$
and let 
$$
\align
y_1&=x_{i_1}^{r_1} x_{i_2}^{r_2}\ldots x_{i_n}^{r_n}\\
y_2&=x_{j_1}^{s_1} x_{j_2}^{s_2}\ldots x_{j_m}^{s_m}
\endalign
$$
be the two positive words involved in the normal form for $x=y_1y_2^{-1}$.
Then the number $N(x)$ of carets for any tree in the reduced tree diagram for 
$x$ is equal to the highest number of carets in the diagrams for $y_1$ and
$y_2$. This number of carets is equal to
$$
\align
N(x)=
\max\,\{&N(y_1),N(y_2)\}\\=\max\,\{&i_k+r_k+r_{k+1}+\ldots+r_n+1,\text{ for }
k=1,2,\ldots,n,\\&j_l+s_l+s_{l+1}+\ldots+s_m+1,\text{ for }l=1,2,\ldots,m\},
\endalign
$$
and it is equivalent to the distance.
\endproclaim

\demo{Proof} The equivalence with the distance proceeds as in Proposition 2,
except the inequalities are now:
$$
\frac{D(x)}4\le N(x)\le D(x)+1.
$$
For the number of carets, we only need to realize that to obtain the diagram 
for $x$ we need to put the two diagrams for $y_1$ and $y_2$ next to each
other with the diagram for $y_2$ reversed, and add carets to the one with
fewer of them. 

The tree diagram obtained will be reduced because of the 
uniqueness of the normal form, again by the results in \cite{8}.\boxx
\enddemo

Note that these results only apply to $F(2)$, where the proofs are
simple and for a positive word in $F(2)$ the number of carets and the $y$-coordinate 
of the last breakpoint coincide. This fact is not true for $F(p)$ if $p>2$, 
but even though those two numbers are different, both are equivalent to 
the distance. For general $F(p)$ we have the following results with analogous
proofs:

\proclaim{Proposition 4} Let 
$$
x=x_{i_1}^{r_1} x_{i_2}^{r_2}\ldots x_{i_n}^{r_n}
$$
be a positive word in $F(p)$. Then the number
$$
N_1(x)=\max\,\left\{\left\lfloor\frac{i_k}{p-1}\right\rfloor
+r_k+r_{k+1}+\ldots+r_n+1,
\text{ for }k=1,2,\ldots,n\right\}
$$
is equal to the number of carets in either tree of the reduced $p$-tree 
diagram for $x$. This number satisfies the inequality
$$
\frac{D(x)}{2(p-1)}\le N_1(x)\le D(x)+1.
$$\boxx
\endproclaim

\proclaim{Proposition 5} Let
$$
x=x_{i_1}^{r_1} x_{i_2}^{r_2}\ldots x_{i_n}^{r_n}
$$
be a positive word in $F(p)$. Then the number
$$
\gather
N_2(x)=\max\,\left\{i_k+r_k(p-1)+r_{k+1}(p-1)+\ldots+r_n(p-1)+1,\right.\\
\left.\text{ for }k=1,2,\ldots,n\right\}
\endgather
$$
is equal to the $y$-coordinate of the last breakpoint of the graph 
for $x$. This number satisfies the inequality
$$
\frac{D(x)}2\le N_2(x)\le D(x)(p-1)+1.
$$\boxx
\endproclaim

And for connection between the number of carets with the word metric 
for a general, not necessarily positive word, we have:

\proclaim{Theorem 6} Let 
$$
x=x_{i_1}^{r_1} x_{i_2}^{r_2}\ldots x_{i_n}^{r_n}
x_{j_m}^{-s_m}\ldots x_{j_2}^{-s_2} x_{j_1}^{-s_1}
$$
be the unique normal form of an element $x\in F(p)$. Then the number
$$
\align
N(x)=\max\,&\left\{\left\lfloor\frac{i_k}{p-1}\right\rfloor
+r_k+r_{k+1}+\ldots+r_n+1,
\text{ for }k=1,2,\ldots,n,\right.\\
&\left.\left\lfloor\frac{j_l}{p-1}\right\rfloor+s_l+s_{l+1}+\ldots+s_m+1,
\text{ for }l=1,2,\ldots,m\right\} 
\endalign
$$
is equal to the number of carets in either tree for the reduced diagram
for $x$. This number is equivalent to the distance. In particular, it 
satisfies the inequality
$$
\frac{D(x)}{4(p-1)}\le N(x)\le D(x)+1.
$$\boxx
\endproclaim

This estimate of the distance in terms of the number of carets in the
trees of the diagram will be used extensively in the next sections, to prove
that several embeddings of a group $F(p)$ in another group $F(q)$ are 
nondistorted.

\heading
3. The embedding of $F(p^k)$ in $F(p)$
\endheading

There are several types of embeddings of groups $F(q)$ as subgroups of groups $F(p)$. The
most natural one is when $q$ is a power of $p$, since then $\ZZ[\frac1{q}]=\ZZ[\frac1p]$. It is
easier to understand these embeddings in terms of carets than in terms of homeomorphisms of
$[0,1]$,  as can be seen in the example $F(4)\subset F(2)$, which we will describe below
in detail. 


The embedding $i:F(4)\rightarrow F(2)$ can be 
understood using the tree diagrams.
Let $T$ be a rooted 4-tree. As seen in section 2, $T$ can be understood
as a subdivision of $[0,1]$ in intervals of length $\frac1{4^k}$. But
clearly, subdividing an interval in four equal parts corresponds to
subdivide the interval first in two parts and then each of these parts
in two more. So given a rooted 4-tree $T$, there is a (unique) rooted 
2-tree $i(T)$ which yields the same subdivision of $[0,1]$. The tree
$i(T)$ can be obtained from $T$ replacing each 4-caret of $T$ by a set of 
three 2-carets in the obvious manner:

\beginpicture
\setcoordinatesystem point at -180 50
\setplotarea x from -180 to 180, y from -40 to 50
\plot -90 0 -60 30 -30 0 /
\plot -70 0 -60 30 -50 0 /
\plot 30 0 40 15 55 30 70 15 80 0 /
\plot 40 15 50 0 /
\plot 60 0 70 15 /
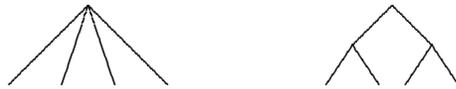
\put {{\bf Fig. 8:} A 4-caret and the set of 2-carets for the embedding of 
$F(4)$ in $F(2)$.} [t] at 0 -10
\endpicture

Now, given an element $x\in F(4)$, represented by the 
reduced tree diagram $(S,T)$,
the element $i(x)\in F(2)$ is represented by the tree diagram $(i(S),i(T))$.
Clearly the diagram $(i(S),i(T))$ will not be reduced in general, so it is
necessary to reduce it. A table with the four generators of $F(4)$
with their corresponding (reduced) images in $F(2)$ is detailed in Figure 9.

Observe that in the process of reducing the diagram $(i(S),i(T))$, some 
2-carets will be eliminated.  But for a set of three 2-carets which corresponds 
to the image of a 4-caret will never be completely erased in the reduction,
because if it were, that would mean that the 4-caret they correspond to
would already have been superfluous. 
So, in the trees $i(S)$ and $i(T)$ there are
at least as many carets as there were in $S$ and in $T$. This provides
the following inequality:
$$
N(x)\le N(i(x)) \le 3N(x).
$$

\beginpicture
\setcoordinatesystem point at -180 360
\setplotarea x from -180 to 180, y from -20 to 360
\plot -145 320 -130 340 -115 320 /
\plot -135 320 -130 340 -125 320 /
\plot -160 300 -145 320 -130 300 /
\plot -150 300 -145 320 -140 300 /
\plot -85 320 -70 340 -55 320 /
\plot -75 320 -70 340 -65 320 /
\plot -70 300 -55 320 -40 300 /
\plot -60 300 -55 320 -50 300 /
\put {$x_0\in F(4)$} [t] at -100 290

\setdots <2pt>
\plot 25 300 30 310 35 300 /
\setsolid
\plot 45 300 50 310 55 300 /
\plot 30 310 40 320 50 310 /
\plot 40 320 45 330 50 320 /
\setdots <2pt>
\plot 60 320 65 330 70 320 /
\setsolid
\plot 45 330 55 340 65 330 /
\setdots <2pt>
\plot 85 320 90 330 95 320 /
\setsolid
\plot 105 320 110 330 115 320 /
\plot 90 330 100 340 110 330 /
\plot 100 300 105 310 110 300 /
\setdots <2pt>
\plot 120 300 125 310 130 300 /
\setsolid
\plot 105 310 115 320 125 310 /
\put {$x_0^2x_1x_2^{-1}\in F(2)$} [t] at 77.5 290

\plot -145 240 -130 260 -115 240 /
\plot -135 240 -130 260 -125 240 /
\plot -150 220 -135 240 -120 220 /
\plot -140 220 -135 240 -130 220 /
\plot -85 240 -70 260 -55 240 /
\plot -75 240 -70 260 -65 240 /
\plot -70 220 -55 240 -40 220 /
\plot -60 220 -55 240 -50 220 /
\put {$x_1\in F(4)$} [t] at -100 210

\plot 35 220 40 230 45 220 /
\setdots <2pt>
\plot 55 220 60 230 65 220 /
\setsolid
\plot 40 230 50 240 60 230 /
\plot 40 240 45 250 50 240 /
\setdots <2pt>
\plot 60 240 65 250 70 240 /
\setsolid
\plot 45 250 55 260 65 250 /
\plot 85 240 90 250 95 240 /
\plot 105 240 110 250 115 240 /
\plot 90 250 100 260 110 250 /
\setdots <2pt>
\plot 100 220 105 230 110 220 /
\plot 120 220 125 230 130 220 /
\setsolid
\plot 105 230 115 240 125 230 /
\put {$x_0x_1^2x_0^{-1}\in F(2)$} [t] at 77.5 210

\plot -145 160 -130 180 -115 160 /
\plot -135 160 -130 180 -125 160 /
\plot -140 140 -125 160 -110 140 /
\plot -130 140 -125 160 -120 140 /
\plot -85 160 -70 180 -55 160 /
\plot -75 160 -70 180 -65 160 /
\plot -70 140 -55 160 -40 140 /
\plot -60 140 -55 160 -50 140 /
\put {$x_2\in F(4)$} [t] at -100 130

\plot 45 140 50 150 55 140 /
\plot 65 140 70 150 75 140 /
\plot 50 150 60 160 70 150 /
\setdots <2pt>
\plot 40 160 45 170 50 160 /
\setsolid
\plot 60 160 65 170 70 160 /
\plot 45 170 55 180 65 170 /
\setdots <2pt>
\plot 85 160 90 170 95 160 /
\setsolid
\plot 105 160 110 170 115 160 /
\plot 90 170 100 180 110 170 /
\plot 100 140 105 150 110 140 /
\plot 120 140 125 150 130 140 /
\plot 105 150 115 160 125 150 /
\put {$x_1^2x_3x_2^{-1}\in F(2)$} [t] at 77.5 130

\plot -145 80 -130 100 -115 80 /
\plot -135 80 -130 100 -125 80 /
\plot -130 60 -115 80 -100 60 /
\plot -120 60 -115 80 -110 60 /
\plot -145 40 -130 60 -115 40 /
\plot -135 40 -130 60 -125 40 /
\plot -85 80 -70 100 -55 80 /
\plot -75 80 -70 100 -65 80 /
\plot -70 60 -55 80 -40 60 /
\plot -60 60 -55 80 -50 60 /
\plot -55 40 -40 60 -25 40 /
\plot -45 40 -40 60 -35 40 /
\put {$x_3\in F(4)$} [t] at -100 30

\plot 55 60 60 70 65 60 /
\setdots <2pt>
\plot 75 60 80 70 85 60 /
\setsolid
\plot 60 70 70 80 80 70 /
\setdots <2pt>
\plot 40 80 45 90 50 80 /
\setsolid
\plot 60 80 65 90 70 80 /
\plot 45 90 55 100 65 90 /
\setdots <2pt>
\plot 40 40 45 50 50 40 /
\setsolid
\plot 60 40 65 50 70 40 /
\plot 45 50 55 60 65 50 /
\setdots <2pt>
\plot 85 80 90 90 95 80 /
\setsolid
\plot 105 80 110 90 115 80 /
\plot 90 90 100 100 110 90 /
\setdots <2pt>
\plot 100 60 105 70 110 60 /
\setsolid
\plot 120 60 125 70 130 60 /
\plot 105 70 115 80 125 70 /
\plot 115 40 120 50 125 40 /
\setdots <2pt>
\plot 135 40 140 50 145 40 /
\setsolid
\plot 120 50 130 60 140 50 /
\put {$x_2^2x_3x_4^{-1}\in F(2)$} [t] at 77.5 30

\put {{\bf Fig. 9:} The images of the generators of $F(4)$ in $F(2)$.}
[t] at 0 10

\endpicture

\noindent By virtue of the results in section 2, 
the number of carets in the reduced
diagram of an element is equivalent to the distance, we obtain that the 
norm of an element in $F(4)$ and of its image in $F(2)$ are equivalent,
so the embedding of $F(4)$ in $F(2)$ is a quasi-isometric embedding.

In the general case, with $i:F(p^k)\longrightarrow F(p)$, the embedding
can also be seen with trees in the exact same way. A $p^k$-caret is now
replaced by 
$$
1+p+p^2+\ldots+p^{k-1}=\frac{p^k-1}{p-1}
$$ 
$p$-carets, and
after the reductions, at least one of the $p$-carets survives per each
$p^k$-caret. This gives the following result:

\proclaim{Theorem 7} The natural embedding of $F(p^k)$ in $F(p)$ is a
quasi-isometric embedding.\boxx
\endproclaim

The inequalities here are
$$
N(x)\le N(i(x))\le\frac{p^k-1}{p-1}N(x).
$$

Note that this embedding is not quasi-onto, since
the image is nowhere near being $\epsilon$-dense. For a general element
in $F(4)$, its image in $F(2)$ will have norm about three times as large.

\heading
4. Embeddings of $F(p)$ in $F(q)$ for any $p$, $q$
\endheading

In \cite{2}, Brin and Guzm\'an proved that any $F(p)$ can be found as a 
subgroup of any $F(q)$. The goal of this section is to use tree diagrams to
understand the behavior of these maps, and prove that all these embeddings are 
quasi-isometric.

To understand how $F(p)$ embeds in $F(q)$, we need to study first a 
particular case. Assume that $p$ and $q$ are such that $q-1$ is a multiple 
of $p-1$. Let $q-1=d(p-1)$. Then we have the following embeddings 
(see \cite{2}):
$$
\gather 
j_1:F(p)\longrightarrow F(q)\\
x_i\mapsto x_{di}
\endgather
$$
and
$$
\gather 
j_2:F(q)\longrightarrow F(p)\\
x_i\mapsto x_i^d.
\endgather
$$
It is easy to verify that these maps preserve the relations of the 
presentations $\Cal P_p$ and $\Cal P_q$. We will give an interpretation 
of these embeddings in terms of tree diagrams, which will be used to 
prove that all these embeddings are quasi-isometric.

To understand the embedding $j_1$, consider first the simplest of these 
embeddings, which embeds $F(2)$ as a subgroup of $F(3)$. 
Observe that the generators 
$x_0$ and $x_1$ of $F(2)$ and $x_0$ and $x_2$ of $F(3)$ have the same number
of carets, and in the same disposition, the only difference is that the
carets for $F(2)$ have two edges, whereas the carets for $F(3)$ have three
(see Figure 10).

So the embedding $j_1$ of $F(2)$ into $F(3)$ can be realized
in terms of the tree diagrams by, given a rooted 2-tree, just
adding a new edge in the middle of every caret to transform it into a 3-caret.
Given an element of $F(2)$ with its tree diagram, we add an edge to every 
caret in both trees to obtain a 3-tree diagram for the image of the given 
element in $F(3)$. And the resulting 3-tree diagram is reduced if and only 
if the starting 2-tree diagram is reduced, since the carets are arranged in
exactly the same pattern. Hence both diagrams have the
same number of carets, and since the number of carets is equivalent to 
the distance, the embedding is a quasi-isometric embedding.

\beginpicture
\setcoordinatesystem point at -180 110
\setplotarea x from -180 to 180, y from -60 to 110
\plot -120 70 -100 90 -90 80 /
\plot -110 80 -100 70 /
\plot -70 80 -60 90 -40 70 /
\plot -60 70 -50 80 /
\put {$x_0\in F(2)$} [t] at -80 60

\plot 30 70 50 90 /
\plot 40 70 40 80 50 70 /
\plot 50 80 50 90 60 80 /
\plot 90 90 110 70 /
\plot 80 80 90 90 90 80 /
\plot 90 70 100 80 100 70 /
\put {$x_0\in F(3)$} [t] at 70 60

\plot -120 20 -110 30 -90 10 /
\plot -120 0 -100 20 /
\plot -110 10 -100 0 /
\plot -70 20 -60 30 -30 0 /
\plot -60 10 -50 20 /
\plot -50 0 -40 10 /
\put {$x_1\in F(2)$} [t] at -80 -10

\plot 40 30 50 20 30 0 /
\plot 30 20 40 30 40 20 /
\plot 50 10 50 20 60 10 /
\plot 40 0 40 10 50 0 /
\plot 90 30 120 0 /
\plot 80 20 90 30 90 20 /
\plot 90 10 100 20 100 10 /
\plot 100 0 110 10 110 0 /
\put {$x_2\in F(3)$} [t] at 70 -10

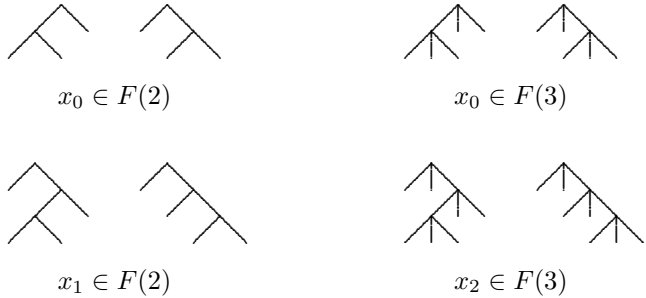
\put {{\bf Fig. 10:} The two generators of $F(2)$ and their images in $F(3)$.}
[t] at 0 -30

\endpicture

It is easy to generalize this interpretation to any $p$, $q$ such that
$q-1=d(p-1)$. Take a $p$-caret
and insert $d-1$ extra edges between every two of the original edges. The
resulting caret has $q$ edges, so by doing this to any rooted $p$-tree
we obtain a rooted $q$-tree with the same number of carets.
For any reduced $p$-tree diagram we obtain the reduced $q$-diagram for the
image of the element it represents, and both diagrams have the same number of
carets. Hence:

\proclaim{Proposition 8} If $q-1=d(p-1)$, the
embedding $j_1$ of $F(p)$ in $F(q)$ is a quasi-isometric embedding.\boxx
\endproclaim

The embedding in the other direction can also be understood in terms of
the tree diagrams. As before, we will study a special case to understand how
the embedding works, and, to make a comparison with the embedding in section
3, we will study the embedding of $F(4)$ into $F(2)$.

To find the image of an element of $F(4)$ in $F(2)$ under the embedding $j_2$, 
we will take the reduced 4-tree diagram for the element, and, as before, 
replace every 4-caret by an arrangement of 2-carets. This time, though, the 
replacement rule will depend on the caret being left, right or interior. The 
replacement rule is as follows: replace every left and interior 4-caret
by a set of three 2-carets where each 2-caret (except the top one) hangs from
the left edge of the caret above; and every right 4-caret by a set of 
three 2-carets where each of the two bottom carets is attached to the 
right edge of the previous one. See figure 11.
The root caret is considered a right caret
for the purposes of this rule. 

\beginpicture
\setcoordinatesystem point at -180 205
\setplotarea x from -180 to 180, y from -5 to 205
\plot -30 155 0 195 30 155 /
\plot -10 155 0 195 10 155 /
\put {1} [t] at -30 150
\put {2} [t] at -10 150
\put {3} [t] at 10 150
\put {4} [t] at 30 150
\plot -75 70 -45 115 -35 100 /
\plot -55 70 -65 85 /
\plot -45 85 -55 100 /
\plot 75 70 45 115 35 100 /
\plot 55 70 65 85 /
\plot 45 85 55 100 /
\put {1} [t] at -75 65
\put {2} [t] at -55 65
\put {3} [t] at -45 80
\put {4} [t] at -35 95
\put {4} [t] at 75 65
\put {3} [t] at 55 65
\put {2} [t] at 45 80
\put {1} [t] at 35 95
\put {left and interior} [t] at -55 45
\put {right} [t] at 55 45
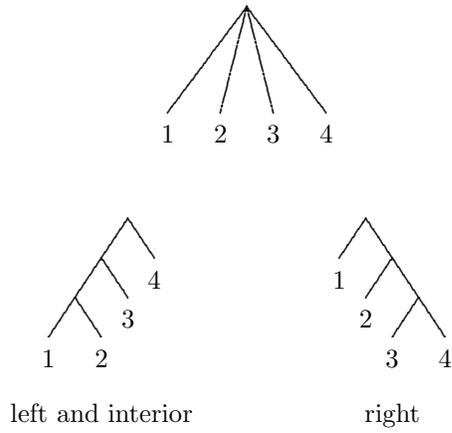
\put {{\bf Fig. 11:} The replacement rule for the embedding $j_2:F(4)
\longrightarrow F(2)$} [t] at 0 15
\endpicture

The motivation for this
choice of replacement rules comes from the connections 
between the rooted tree diagrams and the normal forms of elements, 
stated in \cite{6} and \cite{8}.
Only the left and interior carets contribute to the exponents
of the normal form. Hence the right carets have to be replaced by an 
arrangement made out of right carets to not increase the exponent further.
This caret substitution process will correspond to the desired embedding,
namely $x_1\mapsto x_i^d$.

\beginpicture
\setcoordinatesystem point at -180 440
\setplotarea x from -180 to 180, y from -20 to 440
\plot -145 400 -130 420 -115 400 /
\plot -135 400 -130 420 -125 400 /
\plot -160 380 -145 400 -130 380 /
\plot -150 380 -145 400 -140 380 /
\plot -85 400 -70 420 -55 400 /
\plot -75 400 -70 420 -65 400 /
\plot -70 380 -55 400 -40 380 /
\plot -60 380 -55 400 -50 380 /
\put {$x_0\in F(4)$} [t] at -100 370

\plot 25 380 40 410 50 420 60 410 /
\plot 40 410 45 400 /
\plot 35 400 40 390 /
\plot 30 390 35 380 /
\setdots <2pt>
\plot 55 400 60 410 70 390 /
\plot 60 390 65 400 /
\setsolid
\plot 85 420 90 430 110 390 /
\plot 90 410 95 420 /
\plot 95 400 100 410 /
\plot 100 390 105 400 /
\setdots <2pt>
\plot 105 380 110 390 120 370 /
\plot 110 370 115 380 /
\setsolid
\put {$x_0^3\in F(2)$} [t] at 77.5 360

\plot -145 300 -130 320 -115 300 /
\plot -135 300 -130 320 -125 300 /
\plot -150 280 -135 300 -120 280 /
\plot -140 280 -135 300 -130 280 /
\plot -85 300 -70 320 -55 300 /
\plot -75 300 -70 320 -65 300 /
\plot -70 280 -55 300 -40 280 /
\plot -60 280 -55 300 -50 280 /
\put {$x_1\in F(4)$} [t] at -100 270

\plot 45 315 50 325 55 315 /
\plot 30 275 45 305 55 315 65 305 /
\plot 45 305 50 295 /
\plot 40 295 45 285 /
\plot 35 285 40 275 /
\setdots <2pt>
\plot 60 295 65 305 70 295 /
\setsolid
\plot 85 320 90 330 115 280 /
\plot 90 310 95 320 /
\plot 95 300 100 310 /
\plot 100 290 105 300 /
\plot 105 280 110 290 /
\setdots <2pt>
\plot 110 270 115 280 120 270 /
\setsolid
\put {$x_1^3\in F(2)$} [t] at 77.5 260

\plot -145 200 -130 220 -115 200 /
\plot -135 200 -130 220 -125 200 /
\plot -140 180 -125 200 -110 180 /
\plot -130 180 -125 200 -120 180 /
\plot -85 200 -70 220 -55 200 /
\plot -75 200 -70 220 -65 200 /
\plot -70 180 -55 200 -40 180 /
\plot -60 180 -55 200 -50 180 /
\put {$x_2\in F(4)$} [t] at -100 170

\plot 50 220 55 230 70 200 /
\plot 55 210 60 220 /
\plot 45 170 65 210 /
\plot 60 200 65 190 /
\plot 55 190 60 180 /
\plot 50 180 55 170 /
\plot 85 220 90 230 120 170 /
\plot 90 210 95 220 /
\plot 95 200 100 210 /
\plot 100 190 105 200 /
\plot 105 180 110 190 /
\plot 110 170 115 180 /
\put {$x_2^3\in F(2)$} [t] at 77.5 160

\plot -145 95 -130 115 -115 95 /
\plot -135 95 -130 115 -125 95 /
\plot -130 75 -115 95 -100 75 /
\plot -120 75 -115 95 -110 75 /
\plot -145 55 -130 75 -115 55 /
\plot -135 55 -130 75 -125 55 /
\plot -85 95 -70 115 -55 95 /
\plot -75 95 -70 115 -65 95 /
\plot -70 75 -55 95 -40 75 /
\plot -60 75 -55 95 -50 75 /
\plot -55 55 -40 75 -25 55 /
\plot -45 55 -40 75 -35 55 /
\put {$x_3\in F(4)$} [t] at -100 45

\plot 30 110 35 120 50 90 /
\plot 35 100 40 110 /
\plot 40 90 45 100 /
\plot 25 50 40 80 50 90 60 80 /
\plot 30 60 35 50 /
\plot 35 70 40 60 /
\plot 40 80 45 70 /
\setdots <2pt>
\plot 55 70 60 80 70 60 /
\plot 60 60 65 70 /
\setsolid
\plot 85 120 90 130 125 60 /
\plot 90 110 95 120 /
\plot 95 100 100 110 /
\plot 100 90 105 100 /
\plot 105 80 110 90 /
\plot 110 70 115 80 /
\plot 115 60 120 70 /
\setdots <2pt>
\plot 120 50 125 60 135 40 /
\plot 125 40 130 50 /
\setsolid
\put {$x_3^3\in F(2)$} [t] at 77.5 30

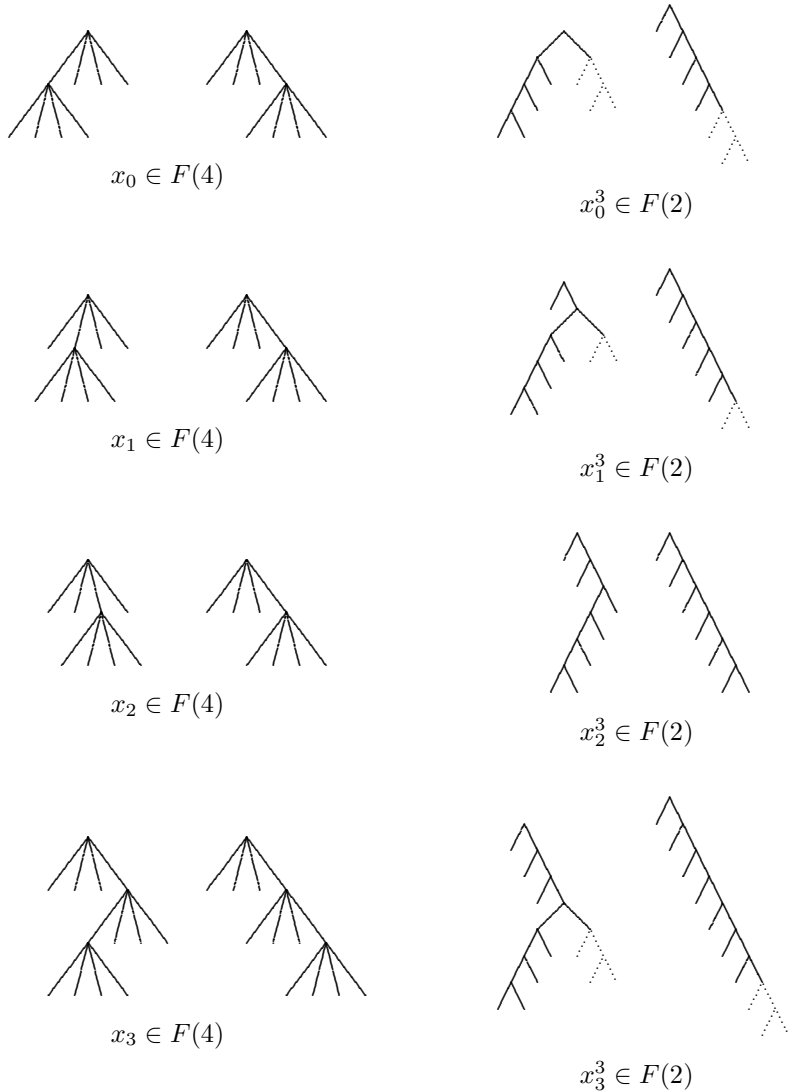
\put {{\bf Fig. 12:} The images of the generators of $F(4)$ in $F(2)$ 
under $j_2$.} [t] at 0 10
\endpicture

Both a 4-caret and either of the two arrangements have four leaves, so the 
arrangements of 2-carets are attached, keeping the root attached
to the corresponding leaf preserving the order. The leaves in figure 11
are numbered from left to right to emphasize this fact.

Observe in figure 12 how this replacement sends the generators $x_i$ of $F(4)$
to the elements $x_i^3$ in $F(2)$. It is easy to check that using the same
replacement rule on a nonreduced diagram yields the same group element. Hence,
the map defined by this replacement rule is a homomorphism, and 
since it corresponds to $j_2$ on the generators, this homomorphism is 
equal to $j_2$. 

The general case is completely analogous to this one. Given $p$ and $q$ with 
$q-1=d(p-1)$, to obtain the embedding $j_2$ from $F(q)$ into $F(p)$ replace
each left or interior $q$-caret by $d$ $p$-carets such that, from the second 
on, each one is attached to its parent on the left leaf; and each right caret
is also replaced by $d$ $p$-carets, all of them on the right side.

It may be neccessary, as it was in the case described in section 3, to reduce
the diagrams obtained after replacing the $q$-carets by the arrangements of
$p$-carets. But observe that, exactly as in the case in section 3, of all 
$p$-carets which originate in the same $q$-caret, at least one survives
after the reduction, because if all were deleted it would mean that the
original $q$-caret would already have been redundant. So as before, the number
of $p$-carets in the image is at most $d$ times the number of $q$-carets 
in the source, and at least the same number. Hence:

\proclaim{Proposition 9} The embedding $j_2$ of $F(q)$ into $F(p)$ is a
quasi-isometric embedding. \boxx
\endproclaim

Combining these two embeddings we can now embed $F(p)$ into $F(q)$, for
any $p,q$.

\proclaim{Theorem 10} For any pair of integers $p,q\ge 2$, the group 
$F(p)$ can be obtained as a subgroup of $F(q)$, and the embedding is
quasi-isometric.
\endproclaim

\demo{Proof} Since $2-1=1$ divides any integer, $F(2)$ embeds in any
$F(p)$, and reciprocally any $F(p)$ embeds in $F(2)$. To embed now
any $F(p)$ in any $F(q)$, embed first $F(p)$ in $F(2)$ and 
then $F(2)$ in $F(q)$.
\boxx
\enddemo

\heading
5. The shift maps
\endheading

One of the most interesting features of the groups $F(p)$ are the shift
maps. The shift map $\phi_p$ for $F(p)$ is defined as the map sending
every generator $x_i$ in the presentation $\Cal P_p$ to $x_{i+1}$.
As we have seen before, $\phi_p$ satisfies
$x_0^{-1}\phi_p(x)x_0=\phi_p^p(x)$. Of special relevance is the map
$\phi_p^{p-1}$, because the relators of $\Cal P_p$ can be written as
$$
x_i^{-1}x_jx_i=\phi_p^{p-1}(x_j),\qquad\text{for }i<j;
$$
so if we take an HNN extension of $F(p)$
by the map $\phi_p^{p-1}$, we obtain another copy of $F(p)$.

The map $\phi_p^{p-1}$ also has significance in the homeomorphisms of
$[0,1]$.
The image of an element $x\in F(p)$ by $\phi_p^{p-1}$ is the identity in
the
interval $\left[0,1-\frac1p\right]$ and has a copy of the graph of $x$
in the interval $\left[1-\frac1p,1\right]$, scaled down by a factor of
$p$.

Also, it is easy to interpret the maps $\phi_p^{p-1}$ in terms of the
rooted $p$-tree diagrams. Given a rooted $p$-tree $T$, consider another
tree $T'$ obtained by taking one single $p$-caret and attaching $T$ by
the root to the rightmost vertex of the caret.

\beginpicture
\setcoordinatesystem point at -180 55
\setplotarea x from -180 to 180, y from -60 to 55

\put {$S$} [b] at -90 30
\put {$T$} [b] at -45 30
\plot -100 15 -90 25 -80 15 /
\plot -100 5 -90 15 -80 5 /
\plot -90 25 -90 5 /
\plot -50 25 -30 5 /
\plot -60 15 -50 25 -50 15 /
\plot -50 5 -40 15 -40 5 /
\put {$x_1\in F(3)$} [t] at -70 -5

\put {$S'$} [b] at 45 35
\put {$T'$} [b] at 100 35
\plot 30 20 40 30 40 20 /
\plot 40 30 60 10 /
\plot 40 10 50 20 50 0 /
\plot 40 0 50 10 60 0 /
\plot 90 30 120 0 /
\plot 80 20 90 30 90 20 /
\plot 90 10 100 20 100 10 /
\plot 100 0 110 10 110 0 /
\put {$\phi_3^2(x_1)=x_3\in F(3)$} [t] at 75 -10

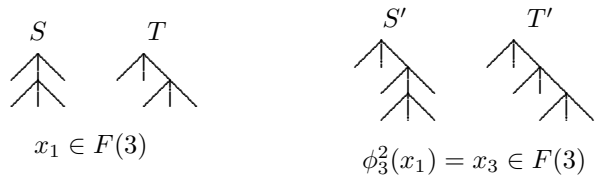
\put {{\bf Fig. 13:} The image of $x_1\in F(3)$ under $\phi_3^2$.} [t] at
10 -30

\endpicture

Using this construction, it is easy to see that 
if $(S,T)$ is the reduced tree diagram for $x\in F(p)$, then
$(S',T')$ is the reduced tree diagram for $\phi_p^{p-1}(x)$.
We can use this interpretation of the shift maps to see that they
behave nicely under the embeddings studied in sections 3 and 4.

\proclaim{Proposition 11} Let
$$
i:F(p^k)\longrightarrow F(p)
$$
be the natural embedding. Then the shift maps satisfy:
$$
i\circ\phi_{p^k}^{p^k-1}=\phi_p^{k(p-1)}\circ i.
$$
\endproclaim
\demo{Proof}
Attaching an extra $p^k$-caret to the top of a rooted tree corresponds
by the embedding $i$ to attaching $\frac{p^k-1}{p-1}$ $p$-carets to the
top of a rooted $p$-tree. But on a reduced $p$-tree diagram, after
attaching these carets to the top of each tree, we see that all the
attached carets except the $k$ rightmost ones will be reduced.

\beginpicture
\setcoordinatesystem point at -180 160
\setplotarea x from -180 to 180, y from -60 to 160

\plot -120 120 -105 140 -90 120 /
\plot -110 120 -105 140 -100 120 /
\plot -125 100 -110 120 -95 100 /
\plot -115 100 -110 120 -105 100 /
\plot -55 120 -40 140 -25 120 /
\plot -45 120 -40 140 -35 120 /
\plot -50 100 -35 120 -20 100 /
\plot -40 100 -35 120 -30 100 /
\put {$x_1 x_2^{-1}\in F(4)$} [t] at -70 90

\plot 40 130 50 140 60 130 /
\plot 35 120 40 130 45 120 /
\plot 55 120 60 130 65 120 /
\plot 35 110 45 120 55 110 /
\plot 30 100 35 110 40 100 /
\plot 50 100 55 110 60 100 /
\plot 105 130 115 140 125 130 /
\plot 100 120 105 130 110 120 /
\plot 120 120 125 130 130 120 /
\plot 120 120 125 130 130 120 /
\plot 110 110 120 120 130 110 /
\plot 105 100 110 110 115 100 /
\plot 125 100 130 110 135 100 /
\put {$i(x_1x_2^{-1})$} [t] at 80 90

\plot -135 40 -120 60 -105 40 /
\plot -125 40 -120 60 -115 40 /
\plot -120 20 -105 40 -90 20 /
\plot -110 20 -105 40 -100 20 /
\plot -125 0 -110 20 -95 0 /
\plot -115 0 -110 20 -105 0 /
\plot -70 40 -55 60 -40 40 /
\plot -60 40 -55 60 -50 40 /
\plot -55 20 -40 40 -25 20 /
\plot -45 20 -40 40 -35 20 /
\plot -50 0 -35 20 -20 0 /
\plot -40 0 -35 20 -30 0 /
\put {$\phi_4^3(x_1x_2^{-1})$} [t] at -80 -10

\plot 25 50 35 60 45 50 /
\setdots <2pt>
\plot 20 40 25 50 30 40 /
\setsolid
\plot 40 40 45 50 50 40 /
\plot 40 30 50 40 60 30 /
\plot 35 20 40 30 45 20 /
\plot 55 20 60 30 65 20 /
\plot 35 10 45 20 55 10 /
\plot 30 0 35 10 40 0 /
\plot 50 0 55 10 60 0 /
\plot 90 50 100 60 110 50 /
\setdots <2pt>
\plot 85 40 90 50 95 40 /
\setsolid
\plot 105 40 110 50 115 40 /
\plot 105 30 115 40 125 30 /
\plot 100 20 105 30 110 20 /
\plot 120 20 125 30 130 20 /
\plot 120 20 125 30 130 20 /
\plot 110 10 120 20 130 10 /
\plot 105 0 110 10 115 0 /
\plot 125 0 130 10 135 0 /
\put {$i(\phi_4^3(x_1x_2^{-1}))=\phi_2^2(i(x_1x_2^{-1}))$}
[t] at 75 -10

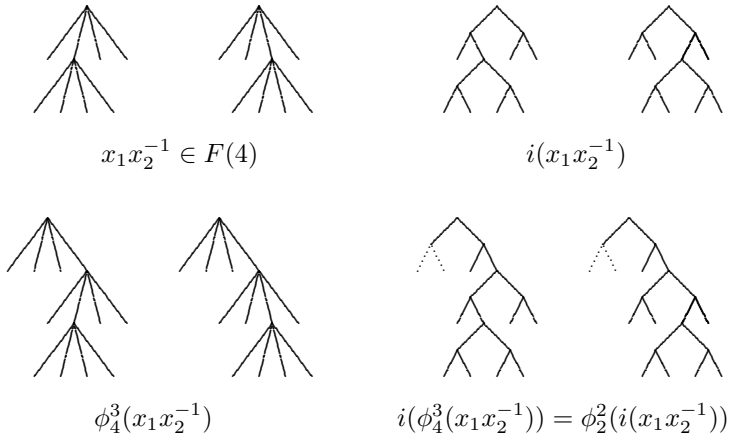
\put {{\bf Fig. 14:} An example of the invariance of the shift maps.}
[t] at 0 -40

\endpicture

So the resulting tree diagram in $F(p)$ is the diagram obtained by
attaching the $k$ carets that are necessary to apply $\phi_p^{k(p-1)}$.
\boxx
\enddemo

\proclaim{Proposition 12} Let
$$
j_1:F(p)\longrightarrow F(q)
$$
be the embedding with $q-1=d(p-1)$. Then the shift maps satisfy:
$$
j_1\circ\phi_{p}^{p-1}=\phi_q^{q-1}\circ j_1.
$$
\endproclaim
\demo{Proof}
Both shift maps act adding a single caret on top of the trees, and this
operation is not affected by adding extra edges to an existing caret,
which is what the embedding $j_1$ demands (see section 4). Clearly the
same result is obtained by adding a $p$-caret and then adding $d-1$
edges between each two than by adding first the sets of $d-1$ edges and
later adding a $q$-caret on top. \boxx
\enddemo

Finally, the shift maps also behave well under the embedding $j_2$:

\proclaim{Proposition 13} Let
$$
j_2:F(q)\longrightarrow F(p)
$$
be the embedding with $q-1=d(p-1)$. Then the shift maps satisfy:
$$
j_2\circ\phi_{q}^{q-1}=\phi_p^{d(p-1)}\circ j_2.
$$
\endproclaim
\demo{Proof} Adding a caret on top of a rooted $q$-tree is what is needed
for $\phi_{q}^{q-1}$. But under the replacement rule 
for $j_2$ this $q$-caret is replaced
by $d$ $p$-carets, all of them right carets, so in $F(p)$ we have applied
the map $\phi_p^{d(p-1)}$.\boxx
\enddemo

\enddocument